\documentclass[preprint,12pt]{elsarticle}

\usepackage{amsmath}
\usepackage{amsfonts}
\usepackage[english]{babel}
\usepackage{graphicx}

\begin{document}

\begin{frontmatter}
\title{Computational implementation of the inverse continuous wavelet transform without a requirement of the admissibility condition}

\author{Eugene B. Postnikov}
\ead{postnicov@gmail.com}
\address{Theoretical Physics Department, Kursk State University, Radishcheva st., 33 Kursk 305000, Russia}

\author{Elena A. Lebedeva}
\ead{ealebedeva2004@gmail.com}
\address{Mathematics and Mechanics Faculty, Saint Petersburg State University,
Universitetsky prospekt, 28, Peterhof,  Saint Petersburg,
 198504, Russia}
\address{Institute of Applied Mathematics and Mechanics, Saint Petersburg State Polytechnic University,
 Polytechnicheskay 29, 195251, Saint Petersburg, Russia}

\author{Anastasia I. Lavrova}
\ead{aurebours@googlemail.com}
\address{Institute of Chemistry and Biology, Immanuel Kant Baltic Federal University, A. Nevskogo str. 14A, Kaliningrad 236041, Russia}

\begin{abstract}
Recently, it has been proven [R. Soc. Open Sci. 1 (2014) 140124] that the continuous wavelet transform with non-admissible kernels (approximate wavelets) allows for an existence of the exact inverse transform. Here we consider the computational possibility for the realization of this approach. We provide modified simpler explanation of the reconstruction formula, restricted on the practical case of real valued finite (or periodic/periodized) samples and the standard (restricted) Morlet wavelet as a practically important example of an approximate wavelet. The provided examples of applications includes the test function and the non-stationary electro-physical signals arising in the problem of neuroscience. 
\end{abstract}


\end{frontmatter}

\section{Introduction}

The continuous wavelet transform (CWT) with the standard (restricted) Morlet wavelet $\psi(\xi)=\exp\left(i\omega_0\xi-\xi^2/2\right)$ is defined as the integral  
\begin{equation}
w(a,b)=\int\limits_{-\infty}^{+\infty}f(t)e^{-i\omega_0\frac{t-b}{a}}
e^{-\frac{(t-b)^2}{2a^2}}\frac{dt}{\sqrt{2\pi a^2}},
\label{CWTM}
\end{equation}
where $a$ and $b$ are the scale and the shift correspondingly, $\omega_0$ is the central frequency.

It is one of the most powerful modern tools of signal processing especially adjusted to the extraction of instant oscillating patterns \cite{MallatBook,AddisonBook} since the transform (\ref{CWTM}) of the harmonic oscillation $f(t)=\exp(i\omega t)$ results in the complex function 
$$
w(a,b)=e^{i\omega b}e^{-\frac{(a\omega-\omega_0)^2}{2}}.
$$
The modulus of this function has a maximum, which allows for the determining of signal's frequency $\omega=\omega_0/a_{max}$ and the phase coincides with signal's one.

The modern applications of the continuous wavelet transform are focused, in particular, on a study of environmental time series \cite{Cazelles2008,Galiana2014}, geo- and astrophysics \cite{Katsavrias2012,Soon2014,Postnikov2009}, biophysics \cite{Meng2013,Worrell2012,Suvichakorn2011} and neuroscience, see an extensive review in the recently published book \cite{WaveletNeuroBook}.

At the same time, the actual problem is not restricted by the search of oscillating patterns localizations. It is important to extract revealed structures from the background consisting a noise, global oscillations and inhomogeneities, etc. \cite{Sheppard2011,Karimi2012,Gu2013,Postnov2014}.

However, the conventional approaches to the inversion of the wavelet transform with the standard Morlet wavelet have some principal difficulties
from the point of view of functional analysis. 
To be applicable in a classical inversion formula, a wavelet function $\psi$ should satisfy the admissibility condition \cite{MallatBook} 
$$
C_{\psi}=\int_{\mathbb{R}}\frac{|\psi^*(\omega)|^2}{|\omega|}d\omega<\infty,
$$  
where the asterisk denotes a complex conjugation.
Then the classical inversion of the wavelet transform $W_{\psi}f=w(a,b)$ with a wavelet function $\psi$ is written as 
$$
f(t)=\frac{1}{C_{\psi}}\int_{\mathbb{R}}\int_{\mathbb{R}} \frac{1}{a}
\psi\left(\frac{t-b}{a}\right) w(a,b) \frac{dadb}{a}.
$$
However, the integral $C_{\psi}$ diverges  for the standard (restricted) Morlet wavelet.

On the other hand, the alternative inversion formula for the CWT 
\begin{equation}
	\frac{1}{\pi}\mathrm{v.p.}\int_{\mathbb{R}}\frac{db}{b-t} \int_{\mathbb{R}}\frac{\partial}{\partial b} w(a,b)\,da = \psi^*(0)\, f(t) 
	\label{rec1} 
\end{equation}
is proven in \cite{Lebedeva2014} under a mild natural conditions on a wavelet function $\psi$.

The principal aim of this paper is to show how this approach can be realized in applications. For this reason, we adapt the proof of (\ref{rec1}) to the case of real-valued functions with a finite support in such a way that its line of reasoning and the result can be straightforwardly used as a practical algorithm for the reconstruction of a function from its wavelet transform with the Morlet wavelet.   
 
The paper is organized as follows. Section 2 shortly discusses the implementation of the CWT with the Morlet wavelet in the application to real-valued functions determined within a finite interval (or periodic on $\mathbb{R}$) and then presents the simple proof the exact reconstruction formula that does not require the wavelet admissibility condition. Its formulation allows for the simple computational realization, which is presented in Section 3. The examples comprise the test one and the solution of the modern practical problem in the field of neuroscience. Neuronal electric activity is characterized by the high non-linearity and complexity, extraction of main oscillatory components and localization of them is very important for reveal of interaction dynamics between different cells in neuronal network \cite{Keef}.
The final section gives some outlooks for the perspectives for the applications of the proposed method in computational physics.

\section{Direct and inverse continuous wavelet transform of periodic functions}

\subsection{CWT expansion}

Let us consider a real-valued function $f(t)$ with the zero mean given by its expansion into the Fourier series 
\begin{equation}
f(t)=\sum\limits_{n=1}^{\infty}A_n\cos(\omega_nt)+B_n\sin(\omega_nt),
\label{fFourier}
\end{equation}
where $\omega_n=2\pi n/T$.

This function is periodic with the period $T$. As well, it can be considered as a function determined within the interval $t\in[0,\,T]$ and periodically extended over all $\mathbb{R}$. The last point of view is applicative in practical computations since one can deal with finite samples only. The assumption of the zero mean is also not restrictive for the problems of computational physics since it can be easily achieved by the exclusion of the averaged value: $f(t)-T^{-1}\int_0^Tf(t)dt\to f(t)$.

For the wavelet analysis of spectral components, it is most convenient to operate with the complex analytical counterpart of the function $f(t)$ obtained via the Hilbert transform $f_a(t)=f(t)+iH\left[f(t)\right]$.

Since $H\left[\cos(x)\right]=\sin(x)$, $H\left[\sin(x)\right]=-\cos(x)$, the Hilbert transform applied to Eq.~(\ref{fFourier}), accompanied with Euler's formula, gives 
\begin{equation}
f_a(t)=\sum\limits_{n=1}^{\infty}\left(A_n-iB_n\right)e^{i\omega_nt}=\sum\limits_{n=1}^{\infty}C_ne^{i\omega_nt}.
\label{anfun}
\end{equation}

Note also that Euler's formula applied directly to Eq.~(\ref{fFourier}) results in the representation
$$
f(t)=\frac{1}{2}\sum\limits_{n=1}^{\infty}C_ne^{i\omega_nt}+C_n^*e^{-i\omega_nt},
$$
where asterisks denotes the complex conjugation. 

Thus, one does not need in practice to apply any special additional procedure for the Hilbert transform of a given sample. It is enough to evaluate the discrete Fourier transform (say, via FFT algorithm) and to cut off elements corresponding to negative frequencies. 

The substitution of the series (\ref{anfun}) into the integral (\ref{CWTM}),
$$
w(a,b)=\sum\limits_{n=1}^{\infty}C_n\int\limits_{-\infty}^{+\infty}e^{i\omega_nt}e^{-i\omega_0\frac{t-b}{a}}
e^{-\frac{(t-b)^2}{2a}}\frac{dt}{\sqrt{2\pi a^2}},
$$
after close form computation of Poisson's integral $\int_{-\infty}^{+\infty}\exp(-x^2)dx=\sqrt{\pi}$, 
results in the desired formula for the continuous wavelet transform with the Morlet wavelet
\begin{equation}
w(a,b)=\sum\limits_{n=1}^{\infty}C_ne^{-\frac{(\omega_na-\omega_0)^2}{2}}e^{i\omega_nb}.
\label{wFourier}
\end{equation}

The expression (\ref{wFourier}) is known and usually applied in the numerical algorithms of CWT computations based on the intermediate FFT, e.g. realised in WaveLab package \cite{WaveLab}. However, the preprocessing cut-off procedure, described above, allows for the shortening of the time/memory consumption because of the twice shortened sample's length in comparison with the standard application of Eq.~(\ref{wFourier}) to the initial function $f(t)$ directly. At the same time, this procedure does not loose any information about $f(t)$ since the complex coefficient $C_n=A_n-iB_n$ contains information on both real-valued ones. 

As well, the representation (\ref{wFourier}) provides the most clear way to obtain the exact inverse continuous wavelet transform.

\subsection{CWT reconstruction}

Now, we aimed to obtain the practically applicable CTW reconstruction formula basing on the series (\ref{wFourier}). Its partial differentiation with respect to the shift $b$ gives the series 
$$
\frac{\partial w(a,b)}{\partial b}=i\sum\limits_{n=1}^{\infty}C_n\omega_ne^{-\frac{(\omega_na-\omega_0)^2}{2}}e^{i\omega_nb}.
$$

Since the included frequencies are strictly positive $\omega_n>0$, the subsequent integral over the non-negative scale half-line leads (using Poisson' integral again) to 
$$
\int\limits_0^{\infty}\frac{\partial w(a,b)}{\partial b}da=i\sqrt{2\pi}\sum\limits_{n=1}^{\infty}
C_ne^{i\omega_nb}=i\sqrt{2\pi}f_a(b).
$$

Therefore, taking into account $C_n=A_n-iB_n$ and applying Euler's formula we get the desired formula for the exact inverse CWT with the Morlet wavelet, which does not depend on the admissibility condition:
\begin{equation}
f(t)=\frac{1}{\sqrt{2\pi}}\mathrm{Im}\left[\int\limits_0^{\infty}\frac{\partial w(a,b)}{\partial b}da\right].
\label{inwCWT}
\end{equation}

Note that it reduces computational complexity in comparison with the full abstract representation (\ref{rec1}) since (\ref{inwCWT}) does not require the additional inverse Hilbert transform. 

Moreover, Eq.~(\ref{inwCWT}) allows for the reconstruction of function's instant features revealed by an exploration of the result of the direct CWT. Since the partial derivative is a local operator, it is possible to consider $w(a,b)$ with $b\in[b_{min}\,b_{max}]$, the interval, which contains the feature, which we are interested in. Thus, restricting (\ref{inwCWT}) on the mentioned $b$, we reconstruct $f(t)$ within $b\in[t_{min}\equiv b_{min}\,t_{max}\equiv b_{max}]$ without boundary disturbances originated from the conventional inversion methods based on convolutions. 

As well, since the Gaussian terms in (\ref{wFourier}) are fast decaying functions with the maximum in $a_n=\omega_0/\omega_n$, the restriction of the integral in (\ref{inwCWT})  on the interval $[a_{min}\,a_{max}]$ reconstruct the oscillating components with the frequencies from $[\omega_0/a_{max}\,\omega_0/a_{min}]$.

Combining this restriction with the choice of some interval of the shift or even considering some non-rectangular region $a=a(b)$, it is possibly to reconstruct and analyse very specific details of a studied function $f(t)$. 

\section{Tests and applications}

\subsection{Test non-stationary oscillations}

As a first numerical example, let us consider the function 
\begin{equation}
f(t)=e^{-4t}\cos(20\pi t)+\chi_{\left[\frac{1}{3},\frac{2}{3}\right]}(t)\sin(40\pi t),
\label{nonstex}
\end{equation}
sampled within the interval $t\in[0,\,1]$. 

Here 
$$
\chi_{[t_b,t_e]}(t)=
\left\{
\begin{array}{lr}
1,&t\in [t_b,\,t_e];\\
0,&t\notin [t_b,\,t_e].
\end{array}
\right.
$$
is the indicator function. Thus, the function (\ref{nonstex}) consists of the decaying oscillations over all studied interval and the stable oscillations of unit amplitude over its middle third only. 

The function (\ref{nonstex}) is sampled in $512=2^9$ equispaced points that allows for the direct application of FFT algorithm in further computations. It should be pointed out that the local boundary values of the function (\ref{nonstex}) are sufficiently different that can result in the drastic boundary effects due to a periodization. Therefore, we used the algorithm, which expand the interval  $[0, \,1]$ into 
$[-0.5, \,1.5]$ 
and function's continuation with boundary reflections \cite{AddisonBook}. The MATLAB code realizing such CWT is presented in Appendix. The central frequency is chosen as 
$\omega_0=2\pi$, 
i.e. maxima lines should be located along $a=0.05$ and $a=0.1$.

\begin{figure}
\includegraphics[width=\textwidth]{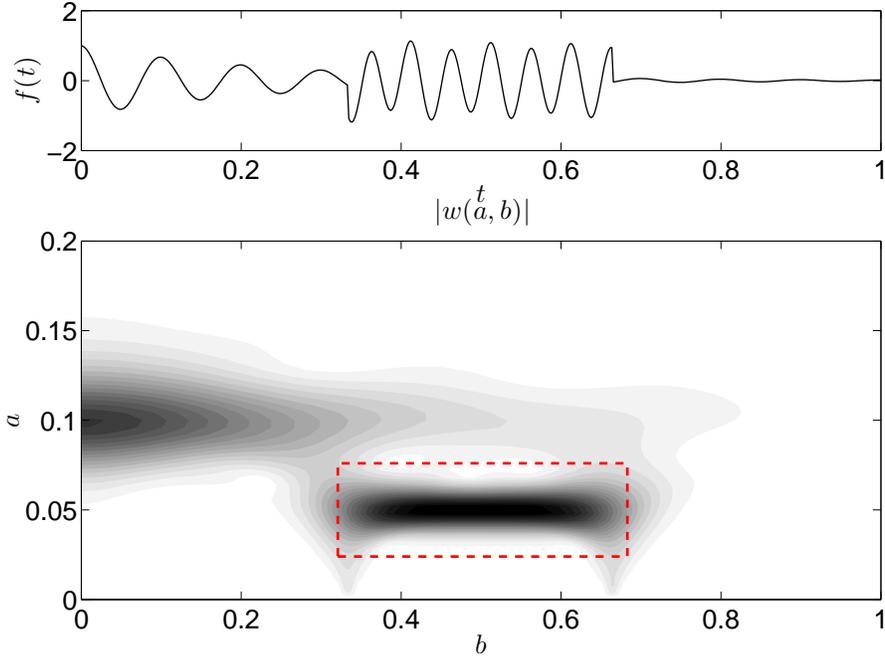}
\caption{The function (\ref{nonstex}) and the modulus of its CWT (dark regions correspond to the larger values of magnitude). The dashed rectangle marks the maximum line corresponding to the localized periodic components with the constant amplitude.}
\label{funccwt}
\end{figure}

\begin{figure}
\includegraphics[width=\textwidth]{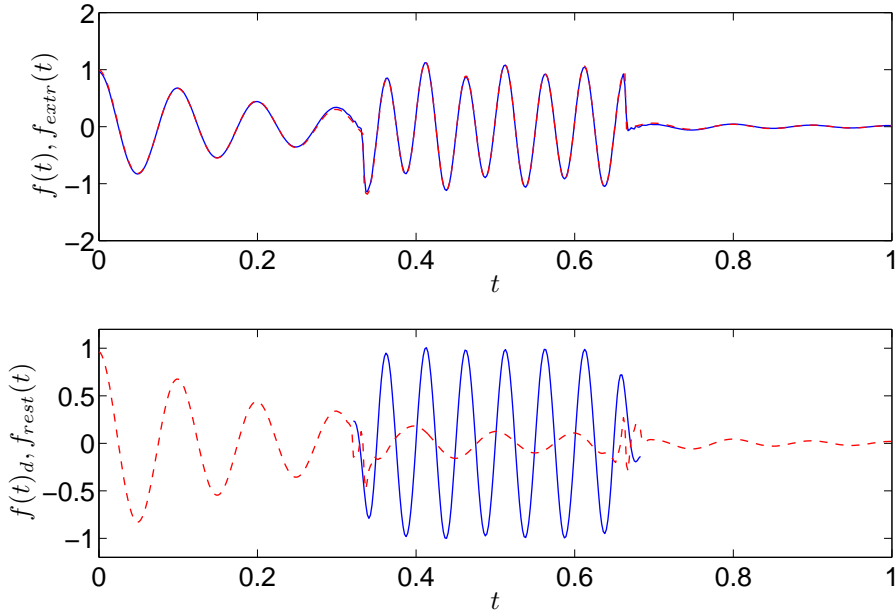}
\caption{The inverse continuous wavelet transform via the formula (\ref{inwCWT}). Upper panel: comparison of the reconstruction (solid line) and the original function (\ref{nonstex}) (dashed line). Lower panel: partial reconstructions of oscillating features of the signal (\ref{nonstex}).}
\label{recfunc}
\end{figure}

Then, the reconstruction formula (\ref{inwCWT}) is applied to this two-dimensional wavelet function $w(a,b)$. The partial derivative with respect to the shift $b$ is realized as the three-point central difference scheme (except the left and the right boundary points, where the forward and  the backward two-point scheme is exploited). The numerical integration is evaluated via the trapezoidal rule. The result of reconstruction for the full shift-scale range is presented in Fig.~\ref{recfunc}. One can see an accurate coincidence of the initial and the reconstructed signals.

As the next step, we process the wavelet transformed function $w(a,b)$ within intervals that bound the localized oscillating features. This should result in the decomposition of the function (\ref{nonstex}) into separate instant oscillatory components. 

As the first step, we extract  $\partial w(a,b)/\partial b$  for the interval $b\in[0.3209,\,0.6829]$ and integrate this part of the obtained function within the limits $a\in[0.0240,\,0.052]$. This region is bounded by the dashed rectangle in Fig.~\ref{funccwt}. The result is drawn by solid line in Fig.~\ref{recfunc}(lower panel). Comparing this localized oscillation with the almost everywhere constant amplitude with the corresponding component of the function (\ref{nonstex}), one can conclude that the proposed methods  reaches its goal. The pattern detected as a localized spot in Fig.~\ref{recfunc}(lower panel) can be reconstructed with a sufficient accuracy using the simple differentiation and integration within the its bounding box. Certainly, it is impossible to avoid boundary disturbances completely, but they are well localized and do not corrupt majority of the reconstructed oscillation. 

The rest component is reconstructed by the subtraction of this local component from the result of full inverse continuous wavelet transform. It is shown as the dashed line in Fig.~\ref{recfunc}(lower panel). One can see that, except the short bursts originated from transient regions, the exponentially decaying oscillation is reconstructed as it should be originally.

\subsection{Neuronal oscillations}

As the application of this method, it would be useful to consider an example of realistic biological oscillations, in particular, the neuronal network dynamics.

We analyse the complex rhythm, which arises prior to epileptic events, so called very fast oscillations (VFO) \cite {Traub}. These oscillations appear spontaneously and can be revealed after the removal of slow baseline fluctuations. Thus, there exists a problem of the reconstruction of various oscillation features from the complex signal representing experimental recordings. Moreover, the considered  method should allows us to determine an approximate time localization for the switching-on and switching-off  of oscillations with different frequencies and their synchronous states. 

\begin{figure}
\includegraphics[width=\textwidth]{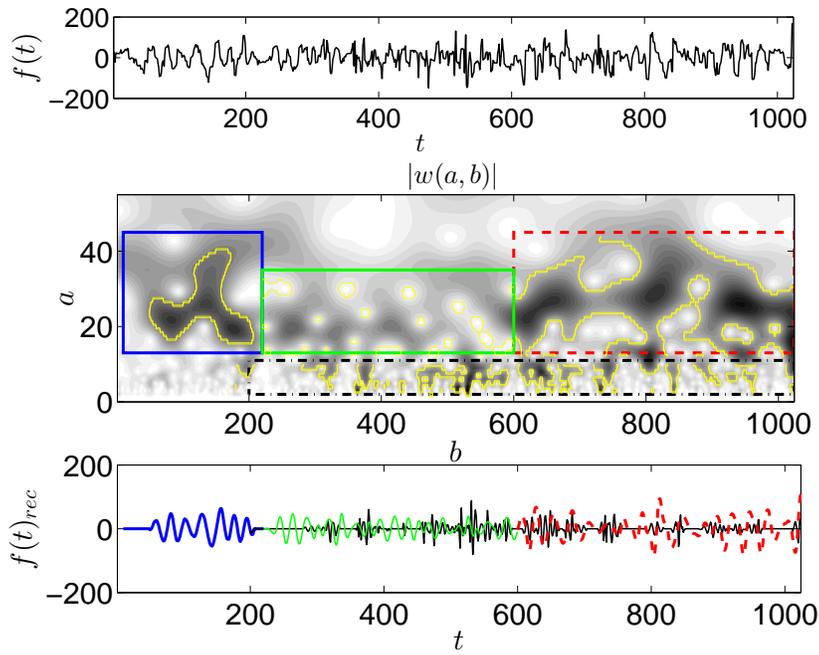}
\caption{The extracellular signal measured in the neocortex layer, its CWT and the partial reconstructions of signal's oscillating features. The plot (upper panel) shows the signal, where VFO occurs (the frequency increases  in the range of time from 250 to 650). Middle panel: Modulus of signal's CWT, where rectangles  mark the regions that contain different oscillating components; the contours bound the exact regions used for the reconstruction. Lower panel: partial reconstructions of the oscillating components contained in the signal. The total time duration corresponds to 4636 ms in the experiment, the maximal amplitude is about 0.14 mV}
\label{ef}
\end{figure}

\begin{figure}
\includegraphics[width=\textwidth]{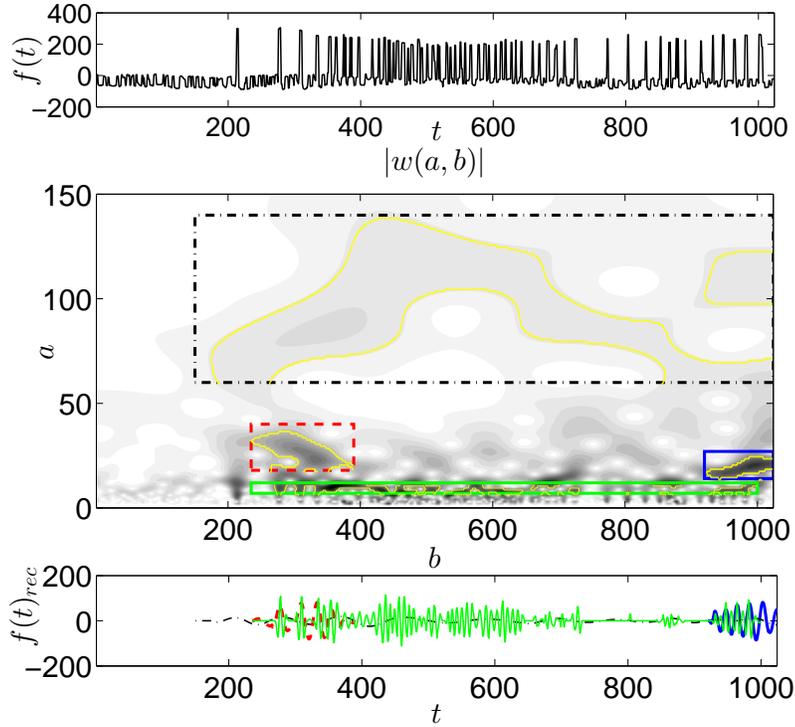}
\caption{The intracellular signal (intrinsic bursting cell \cite {Traub}), its CWT and the partial reconstructions of intracellular signal's oscillating features. The plot (upper panel) shows the signal, where VFO occurs (the frequency increases  in the range of time from 220 to 650). Middle panel: Modulus of signal's CWT, where rectangles  mark the regions that contain different oscillating components; the contours bound the exact regions used for the reconstruction. Lower panel: partial reconstructions of the oscillating components contained in the signal. The total time duration corresponds to 4636 ms in the experiment, the maximal amplitude is about 90 mV}

\label{1cell}
\end{figure}

The experimental data are taken from the work \cite{Traub}. The potential oscillations  are measured in   rat's neocortex: i) the extracellular recordings from the neocortex layer (see Fig.~\ref{ef}, upper panel)  and ii) the intracellular recordings from the pyramidal cell, so-called intrinsically bursting cell, as shown in  Fig.~\ref{1cell}, upper panel. All numerical values are presented in dimensionless units, which can be rescaled to the experimental values of the time and the potential (Fig.1 in \cite{Traub}) as follows:  1 time unit in Figs.~\ref{1cell},\ref{ef} corresponds  to  4.5~ms; 1 voltage unit corresponds to 0.25~mV in Fig.~\ref{1cell} and to 0.5~$\mu$V in Fig.~\ref{ef}.

The visual exploration of the wavelet modulus plots, Figs~\ref{1cell},\ref{ef} (middle panels), allows for concluding that there are at least four areas where oscillations frequencies are different. On the other hand, the corresponding maxima have not such a regular shape as in the the idealized test example considered above. For this reason, we introduce the procedure, which extracts the irregular regions as 
$$
f(t)=\frac{1}{\sqrt{2\pi}}\mathrm{Im}\left[\int\limits_0^{\infty}C(a,b)\frac{\partial w(a,b)}{\partial b}da\right],
$$
where
$$
C(a,b)=\theta\left(|w(a,b)|-L\cdot\mathrm{max}(|w(a,b)|)\right)
$$
is the mask with the cut-off threshold $L\in[0,\,1]$. Here 
$$
\theta(\xi)=\left\{
\begin{array}{lr}
1,\,\xi>0\\
0,\,elsewhere.
\end{array}
\right.
$$
is the Heaviside function. 

Thus, the reconstructed features are the results of the inverse continuous wavelet transform applied to the regions bounded to be the contours defined by the equality $|w(a,b)|=L\cdot\mathrm{max}(|w(a,b)|$.

This reconstruction within contours allows for discussing the types of dynamics, which are characterized by different frequencies (see Figs~\ref{1cell},\ref{ef}, lower panels). In the extracellular signal, as shown in  Fig~\ref{ef}, there are oscillations with two main frequencies: approximately 6 Hz (within the region marked by solid and dashed lines, red, blue and green in the colour online version) and 30 Hz (within the black dash-dotted contour). It should be noted, that the ``fast signal'' (thin black line in the lower panel in Figs~\ref{ef}) is not so strong and arises spontaneously at the time about 250 reaching the maximum of amplitude in the range of time from 450 to 650. On the other hand, the ``fast dynamics'' within the cell signal (thin line (green in the colour online version), Fig~\ref{1cell}, lower panel) starts earlier at time about 220 and has a higher amplitude. It breaks abruptly and renews at the time $t$=940. 

Note that our method, which explicitly extracts the main frequency components within their regions of switch on and off allows for avoiding their mixture with the subthreshold oscillations (so called spikelets), which results in the averaged frequency of ``fast'' signal determined as $118\pm 10$ Hz in \cite {Traub}.

In comparison with the extracellular signal, the  ``slow'' dynamics corresponds to the oscillations of 2 Hz (black dash-dotted line in Fig~\ref{1cell}, lower panel) with a very small amplitude, which do not have breaks within the whole time range of time, to the oscillations of 6 Hz (thick solid (red and blue in the colour online version) lines in Fig~\ref{1cell}, lower panel)). These oscillations exist within the  time intervals from 250 to 400 and from 940 to 1024 only.

\section{Conclusion and outlook}

    In this paper, we have proposed simple computational procedure for the inverse continuous wavelet transform that allows for processing of oscillating signals, e.g. the extraction of main frequency components. The modern methods of physical research, especially connected with electrophysiology, provides data of high complexity, and, correspondingly, their processing requires a development of new approached deeply based of methods of functional analysis \cite{Siddiqi2003book}. The wavelet methods play an important role among them. 
At the same time, they need to be adapted to needs and perceptions of computational community. 		

Thus, we propose new derivation of the reconstruction formula, which, in contrast to the theorem in \cite{Lebedeva2014} provides the direct way to a computational realization and applications. 

As an example, we have considered the voltage signals measured in neuronal system. 
 A neuronal activity is characterized by the both non-linear periodicity and complexity in time and space \cite{Keef, Traub,Buz,Cun}. Firing patterns could be subdivided into the  different frequency bands \cite{Buz,Traub,Cun} and connected  with the spatial periodicity \cite{Keef} that provides an important information on different types of neurons at pathological processes and a normal state. 
    
We have considered the application of this approach to the real experimental data \cite{Traub}, which describe the spontaneous emergence of very fast oscillations at a seizure.  The approach considered in the present work allows for extracting main periodical components from two signals (extracellular and intracellular). The time moments of ``switching on'' of fast oscillations is obtained as well. This information  could be useful for the understanding of a network structure, in particular, a number of oscillating elements and the kind of interconnection between them.

However, the results are not restricted by the problems of electrophysical physiology since the proposed approach and its mathematical and computational background  provide an opportunity to be applied to a much more wider class of complex oscillating patterns. 

\vspace{-5mm}
\section*{Acknowledgment}

EBP is partially supported by grant no. 1391 of the Ministry of Education and Science of the Russian Federation within the basic part of research funding no. 2014/349 assigned to Kursk State University.
EAL is supported by the Russian Foundation for Basic Research, grant No. 15-01-05796 and by the grant No. 9.38.198.2015 of Saint Petersburg State University.
AIL is supported by grant no. 14.575.21.0073, code RFMEFI57514X0073 of the Ministry of Education and Science of the Russian Federation

%\bibliographystyle{elsarticle-num}
%\bibliography{bibInvCWT}

\section*{Appendix}

MATLAB code, which generates the example represented in Figures~\ref{funccwt},~\ref{recfunc}.

\begin{verbatim}
% Determination of the function
N=512;
t=linspace(0,1,N);
f1=cos(20*pi*t);
f2=exp(-4*t);
f3=sin(40*pi*t).*ramp(t,t(round(N/3)),t(round(2*N/3)));
f=f1.*f2+f3;
% Boundary continuation (extension)
[te,fe]=fcontin(t,f);
%Tranform to analytical function
fa=sAnalytic(fe);
\end{verbatim} 

The corresponding functions, which realize extended equispaced samples via the boundary continuation
\begin{verbatim}
function [tp,fp]=fcontin(t,f);

n=length(t);tv(:,1)=t;fv(:,1)=f;
%Interpolation into power-two sequences:
N=2^(ceil(log2(n)));
ti(:,1)=linspace(tv(1),tv(end),N);
fi=interp1q(tv,fv,ti);
%Boundary continuation
N=2^(ceil(log2(n)));
dt=ti(2)-ti(1);
tp=zeros(2*N,1);
tp(1:N/2)=linspace(ti(1)-N*dt/2,-dt,N/2);
tp(N/2+1:1.5*N)=ti;
tp(1.5*N+1:2*N)=linspace(ti(end)+dt,ti(end)+N*dt/2,N/2);
fi=interp1q(tv,fv,ti);
% Boundary reflection
fp(1:N/2)=flipud(conj(fi(2:N/2+1)));
fp(N/2+1:1.5*N)=fi;
fp(1.5*N+1:2*N)=flipud(conj(fi(N/2:N-1)));
\end{verbatim}
and by the cut off of negative frequencies that forms an analytic function: 
\begin{verbatim} 
function fa=sAnalytic(f);

N=length(f);
%Cut-off of negative frequencies
F=fft(f);
F=[2*F(1:N/2),zeros(1,N/2)];
%Output: the analytic signal
fa=ifft(F);
\end{verbatim}

MATLAB function for the CWT with the Morlet wavelet in the amplitude norm, which should be applied to the formed extended sample
\begin{verbatim}
function w=fftMorlet(t,fp,a,omega0);

N=length(t);
%Fourier transform 
F=fft(fp);
nrm=2*pi/(t(end)-t(1));
omega_=([(0:(N/2)) (((-N/2)+1):-1)])*nrm;
%Convolution
if a(1)==0
    w(1,:)=fp*exp(-omega0^2/2);
    k1=2;
else
    k1=1;
end
for k=k1:length(a);
    omega_s=a(k)*omega_;
    window=exp(-(omega_s-omega0).^2/2);
    cnv(k,:)=window.*F;
    w(k,:)=ifft(cnv(k,:));
end
\end{verbatim}

Fig.~\ref{funccwt}(lower panel) presents \verb!abs(w)! is obtained as
\begin{verbatim}
omega0=2*pi;
a=linspace(0,0.2,51);
we=fftMorlet(te,fa,a,omega0);
ti=te(0.5*N+1:1.5*N);
w=we(:,0.5*N+1:1.5*N);
\end{verbatim}

The reconstruction procedure is realized in MATLAB as the function
\begin{verbatim}
function iw=invMorlet(t,a,w);

N=length(t);
d2t=t(3)-t(1);
dw(:,2:N-1)=(w(:,3:end)-w(:,1:end-2))/d2t;
dw(:,1)=(w(:,2)-w(:,1))/(t(2)-t(1));
dw(:,N)=(w(:,N)-w(:,N-1))/(t(N)-t(N-1));
iw=imag(trapz(a,dw))/sqrt(2*pi);
\end{verbatim}

\end{document}